\documentclass[12pt]{amsart}
\usepackage{pifont}
\usepackage{mathrsfs}
\usepackage{geometry}
\usepackage{titletoc}
\usepackage{stix2}
\usepackage{amsmath}
\usepackage{amssymb} %to produce blackboard bold characters(\supsetneq)
\usepackage{enumitem} %lay­out of the three ba­sic list en­vi­ron­ments
\usepackage{mathtools} %an extension package to amsmath + loading amsmath
\usepackage[table]{xcolor} %color of tables
\usepackage[all]{xy} %xymatrix commutative diagram
\usepackage{tikz} %pictures
\usepackage{tikz-cd}%implies
\usepackage{indentfirst} %Indent first paragraph
\usepackage{babel} %lan­guages
\usepackage{setspace} %set­ting the spac­ing be­tween lines

\usepackage[colorlinks,linkcolor=red,anchorcolor=green,citecolor=blue]{hyperref} %color of hyperlinks
\hypersetup{linktocpage = true} %color of hyperlinks of TOC

\usepackage{rotating} %Rotate an image, table or paragraph
\usepackage{ytableau} %young tableau
\usepackage{longtable} %long table
\newcolumntype{M}[1]{>{\centering\arraybackslash}m{#1}} %define dimension for long stable

\usepackage{fancyhdr}

\newcommand\cH{{\mathcal H}}

\newcommand\cV{{\mathcal V}}

\newcommand\bp{{\bar\partial}}

\theoremstyle{plain}
\newtheorem{thm}{Theorem}[section]
\newtheorem{lemma}[thm]{Lemma}
\newtheorem{prop}[thm]{Proposition}
\newtheorem{cor}[thm]{Corollary}
\newtheorem{defn}[thm]{Definition}

\newtheorem{problem}[thm]{Problem}

\theoremstyle{definition}
\newtheorem{example}[thm]{Example}
\newtheorem{remark}[thm]{Remark}

\newcommand{\btheorem}{\begin{thm}}
	\newcommand{\etheorem}{\end{thm}}
\newcommand{\bproposition}{\begin{prop}}
	\newcommand{\eproposition}{\end{prop}}
\newcommand{\bdefinition}{\begin{defn}}
	\newcommand{\edefinition}{\end{defn}}
\newcommand{\bcorollary}{\begin{cor}}
	\newcommand{\ecorollary}{\end{cor}}
\newcommand{\bproof}{\begin{proof}}
	\newcommand{\eproof}{\end{proof}}
\newcommand{\bremark}{\begin{remark}}
	\newcommand{\eremark}{\end{remark}}
\newcommand{\eexample}{\end{example}}
\newcommand{\bexample}{\begin{example}}

\newcommand{\elemma}{\end{lemma}}
\newcommand{\blemma}{\begin{lemma}}

\newcommand{\la}{\langle}
\newcommand{\ra}{\rangle}
\newcommand{\sq}{\sqrt{-1}}

\newcommand{\p}{\partial}

\renewcommand{\bar}{\overline}
\newcommand{\eps}{\varepsilon}

\renewcommand{\phi}{\varphi}

\newcommand{\beq}{\begin{equation}}
\newcommand{\eeq}{\end{equation}}
\newcommand{\ee}{\end{eqnarray*}}
\newcommand{\be}{\begin{eqnarray*}}

\newcommand{\bd}{\begin{enumerate}}
	\newcommand{\ed}{\end{enumerate}}

\renewcommand{\hat}{\widehat}
\renewcommand{\tilde}{\widetilde}

\newcommand{\qtq}[1]{\quad\mbox{#1}\quad}
\renewcommand{\bp}{\bar{\partial}}

\renewcommand{\S}{{\mathbb S}}

\newcommand{\C}{{\mathbb C}}

\renewcommand{\P}{{\mathbb P}}

\newcommand{\R}{{\mathbb R}}

\newcommand{\LL}{\left\langle}
\newcommand{\RL}{\right\rangle}

\newcommand{\om}{\omega}

\setlist[itemize]{leftmargin=*}
\setlist[enumerate]{leftmargin=*}

\numberwithin{equation}{section} %numbering of equations
\makeatletter

\usepackage{fancyhdr}
\title{Conjugate radius, volume comparison and rigidity}
\author{Zhiyao Xiong}
\author{Xiaokui Yang}

\address{Zhiyao Xiong, Department of Mathematics, Tsinghua University, Beijing, 100084, China}
\email{xiongzy22@mails.tsinghua.edu.cn}

\address{Xiaokui Yang, Department of Mathematics and Yau Mathematical Sciences Center, Tsinghua University, Beijing, 100084, China}
\email{xkyang@mail.tsinghua.edu.cn}

\begin{document}
	
	\begin{abstract} 
		In this paper, we prove conjugate radius estimate, volume comparison and rigidity theorems   for  K\"ahler manifolds with various curvature conditions.  
	\end{abstract}

	\maketitle

	\section{Introduction}
		
	Comparison theorems are crucial tools for understanding geometric concepts in differential geometry. 
	Let $(M,g)$ be a complete $n$-dimensional Riemannian manifold with Ricci curvature $\mathrm{Ric}(g)\geq (n-1)g$.
	Myers \cite{Mye41} established the diameter comparison theorem that $\mathrm{diam}(M,g)\leq \mathrm{diam}(\S^n,g_{\mathrm{can}})=\pi$.  Moreover,
Cheng \cite{Che75} obtained the diameter rigidity theorem, 
		which states that if the diameter $\mathrm{diam}(M,g)=\mathrm{diam}(\S^n,g_{\mathrm{can}})$, then $(M,g)$ is isometric to the round sphere. 
	Furthermore, the Bishop-Gromov volume comparison theorem (e.g. \cite{BC64}, \cite{Gro07}, \cite{CE08}) asserts that $\mathrm{Vol}(M,g)\leq \mathrm{Vol}(\S^n,g_{\mathrm{can}})$,
		and the identity holds if and only if $(M,g)$ is isometric to the round sphere. 
	In \cite{CC97}, Cheeger and Colding obtained similar rigidity theorems for volume gaps.
	For more details along this comprehensive topic, we refer to \cite{Wei07} and the references therein.\\

	There are many notable extensions on complete K\"ahler manifolds. 
	For instance, Li and Wang \cite{LW05} obtained diameter comparison and volume comparison theorems 
		in the case that the holomorphic bisectional curvature has a positive lower bound $\mathrm{HBSC}\geq 1$. 
	More recently, Datar and Seshadri \cite{DS23} established the diameter rigidity theorem, 
		which states that if $\mathrm{HBSC}\geq 1$ and $\mathrm{diam}(M,g)=\mathrm{diam}(\C\P^n, \omega_{\mathrm{FS}})$, 
		then $(M,\omega_g)$ is isometrically biholomorphic to $(\C\P^n,\omega_{\mathrm{FS}})$. 
	This is achieved by using Siu-Yau's solution to the Frankel conjecture \cite{SY80} and an interesting monotonicity formula for Lelong numbers on $\C\P^n$ (\cite{Lot21}).
	Similar results were proved in \cite{TY12} and \cite{LY18} with some extra conditions. 
	On the other hand, utilizing entirely different techniques from algebraic geometry (e.g. \cite{Fuj18}),
		Zhang \cite{Zha22} obtained volume comparison and rigidity theorems under the assumption $\mathrm{Ric}(\om)\ge (n+1)\om$.\\

	It is also an interesting topic to investigate diameter comparsion, volume comparison and rigidity theorems  for complete K\"ahler manifolds with positive holomorphic sectional curvature.
	Tsukamoto  proved in \cite{Tsu57} that if a complete K\"ahler manifold $(M,\omega_g)$ has holomorphic sectional curvature $\mathrm{HSC}\ge 2$, 
		then $M$ is compact, simply connected, and $\mathrm{diam}(M,g)\le \mathrm{diam}(\C\P^n,\om_{\mathrm{FS}})$.
	Recently, Ni and Zheng \cite{NZ18} obtained interesting Laplacian comparison and volume comparison theorems by assuming $\mathrm{HSC}\ge 2$ and orthogonal Ricci curvature $\mathrm{Ric}^\perp\ge (n-1)$. \\		
		
	In this paper, we derive volume comparison and rigidity theorems for Kähler manifolds under various curvature conditions. 
	Additionally, we establish  conjugate radius and injectivity radius estimates and the corresponding rigidity theorems.\\
			
  For the reader's convenience, we fix some terminologies.	Let $(M,g)$ be a complete Riemannian manifold. For a unit vector $v \in T_pM$, 
	$c_v$ is the smallest number $t_0>0$ such that $\gamma(t_0)$ is conjugate to $p$ along the geodesic $\gamma(t)=\exp_p(tv)$.  
	The {conjugate radius of $p$} and the {conjugate radius of $(M,g)$} are defined as 
	\[ 
	\mathrm{conj}(M,p):=\inf_{v\in T_pM,\,|v|=1}c_v\qtq{and}\mathrm{conj}(M,g):=\inf_{p\in M}\mathrm{conj}(M,p).
	\] 
	The first main result of this paper is the following volume comparison and rigidity theorem for K\"ahler manifolds with positive holomorphic sectional curvature.
	\btheorem  \label{main4} 
	Let $(M,\omega_g)$ be a complete K\"ahler manifold with $\mathrm{HSC}\ge 2$. If there exists some point $p\in M$ such that  $\mathrm{conj}(M,p)\geq \pi/\sqrt{2}$, then 
	\beq 
	\mathrm{Vol}\left(M,\omega_g\right)\le \mathrm{Vol}\left(\C\P^n, \omega_{\mathrm{FS}}\right), \label{vol comparison}
	\eeq  
	and the identity holds if and only if $(M,\omega_g)$ is isometrically biholomorphic to $(\C\P^n,\omega_{\mathrm{FS}})$.
	\etheorem 
	
	\noindent This result is obtained by utilizing relationships between the RC-positivity proposed in \cite{Yan18} 
		and the conjugate radius estimate derived from the index theorem. The second named author established in \cite{Yan18} that compact K\"ahler manifolds with positive holomorphic sectional curvature are projective and rationally connected, which confirmed affirmatively a conjecture proposed by S.-T. Yau in \cite[Problem~47]{Yau82}, and such manifolds are not necessarily $\C\P^n$.
	The main difficulty in achieving volume comparison and rigidity theorems for compact K\"ahler manifolds with $\mathrm{HSC}>0$ is that the holomorphic sectional curvature is too weak to obtain Laplacian comparison type theorems (see Problem \ref{HSC}). Actually, we derive extra curvature relation from the lower bound of the conjugate radius  at some point. 
	Moreover, we establish (global) conjugate radius and injective radius estimates for such manifolds.
	\btheorem  \label{main5} 
	Let $(M,\omega_g)$ be a complete K\"ahler manifold with $\mathrm{HSC}\ge 2$. 
	Then \beq \mathrm{inj}(M,g)=\mathrm{conj}(M,g)\leq \frac{\pi}{\sqrt{2}},\label{inj=conj}\eeq  
	and the identity holds if and only if $(M,\omega_g)$ is isometrically biholomorphic to $(\C\P^n,\omega_{\mathrm{FS}})$.
	\etheorem
	
	\noindent By using perturbations of $(\C\P^n, \omega_{\mathrm{FS}})$, 
	it is easy to see that there exists a compact K\"ahler manifold $M$ with $\mathrm{HSC}\ge 2$ and $\mathrm{conj}(M,g)\leq \eps<\pi/\sqrt{2}$, 
	but there exists a point $p\in M$ such that $\mathrm{conj}(M,p)\geq \pi/\sqrt{2}$.
	Hence, the conditions in Theorem \ref{main4} cannot be implied by those in Theorem \ref{main5}. 
	On the other hand, Theorem \ref{main5} is a generalization of classical results obtained in \cite{Kli59} and \cite{Gre63} (see also \cite{AM94}) for Riemannian manifolds. 
	Indeed, if $(M,g)$ is a compact Riemannian manifold with scalar curvature $\geq n(n-1)$, 
	Green \cite{Gre63} proved that the conjugate radius $\mathrm{conj}(M,g)\le \pi$, and the identity holds if and only if $(M,g)$ is isometric to the round sphere. 
	We also obtain the following extension in K\"ahler geometry.
	
	\btheorem\label{main1}
	Let $(M,\omega_g)$ be a compact K\"ahler manifold  of complex dimension $n$. Then 
	\beq 
	\frac{4a^2}{\pi (n+1)!}\int_M c_1(M)\wedge [\omega_g]^{n-1} \le \mathrm{Vol}(M,\omega_g),\label{formula key conclusion}
	\eeq 
	where $a$ is the conjugate radius of $(M,g)$.  Moreover, the identity holds if and only if $(M,\omega_g)$ is isometrically biholomorphic to $\left(\C\P^n,\frac{2a^2}{\pi^2}\omega_{\mathrm{FS}}\right)$.
	\etheorem	
	\noindent It is well-known that if $(M,g)$ is a complete Riemannian manifold with non-positive sectional curvature, then $\mathrm{conj}(M,g)=+\infty$.  
	In Theorem \ref{main1}, the conjugate radius $a$ can also be $+\infty$, and in this case we have $\int_M c_1(M)\wedge [\omega_g]^{n-1} \leq 0$. 
	Note also that Zhu established in \cite{Zhu22} some interesting results on the geometry of positive scalar curvature on complete non-compact Riemannian manifolds with non-negative Ricci curvature, 
		which can also be extended to K\"ahler manifolds by using (total) scalar curvature.
	As an application of Theorem \ref{main1}, one has 
	\bcorollary \label{main2}
	Let $(M,\omega_g)$ be a compact K\"ahler manifold of complex dimension $n$.
	If the scalar curvature of $(M,\omega_g)$  satisfies $s\ge n(n+1)$, then 
	\beq  \mathrm{conj}(M,g)\le \frac{\pi}{\sqrt{2}},\eeq
	and the identity holds if and only if $(M,\omega_g)$ is isometrically biholomorphic to $(\C\P^n,\omega_{\mathrm{FS}})$.
	\ecorollary	
	
	\noindent 
	As another application of Theorem \ref{main1}, we give a criterion for finiteness of $\mathrm{conj}(M,g)$ by using RC-positivity.
		
	\btheorem  \label{main3}
	Let $M$ be a compact K\"ahler manifold. 
	If the anti-canonical line bundle $K_M^{-1}$ is RC-positive,
	then for any K\"ahler metric $\omega_g$ on $M$,  $$\mathrm{conj}(M,g)<+\infty.$$
	\etheorem 

	\noindent Recall that the anti-canonical line bundle $K^{-1}_M$ of a compact complex manifold $M$ is called RC-positive 
		if there exists a Hermitian metric $\omega$ on $TM$ such that the first Chern-Ricci curvature $\mathrm{Ric}^{(1)}(\omega)$ has a positive eigenvalue at each point $p\in M$. 
	It is proved in \cite{Yan19a} and \cite{Yan19b} that $K_M^{-1}$ is RC-positive if and only if $K_M$ is not a pseudo-effective line bundle. 
	Consequently, any K\"ahler metric on a uniruled algebraic manifold has finite conjugate radius. 
	We also observe that the converse of Theorem \ref{main3} is not valid in general. 
	Actually, if $M$ is a complete intersection of two generic hypersurfaces in $\C\P^4$ whose degrees are greater than $35$, 
	it is shown in \cite[Corollary~4.13]{Bro14} that $M$ has ample cotangent bundle and so $K_M$ is pseudo-effective.
	Moreover, since $M$ is simply connected (\cite[pp.~221--222]{Sha13}), any metric $g$ on $M$ must have finite conjugate radius. 
	Otherwise, $M$ would be diffeomorphic to $\R^4$.\\

\	Finally, we  establish volume comparison and rigidity theorems for complete K\"ahler manifolds with positive orthogonal holomorphic bisectional curvature (OHBSC), 
	which generalize  results in \cite{LW05}.

	\btheorem\label{OHBSC}	Let $(M,\om_g)$ be a complete K\"ahler manifold of dimension $n\ge 2$. If  $\mathrm{OHBSC}\geq 1$, then $M$ is compact and
	\beq 
	\mathrm{Vol}(M,\om_g)\le  \mathrm{Vol}(\C\P^n,\om_{\mathrm{FS}}),
	\label{formula HBSC vol upper}
	\eeq
	and the identity holds if and only if $(M,\om_g)$ is isometrically biholomorphic to $(\C\P^n,\om_{\mathrm{FS}})$.
	\etheorem 
	
	\noindent	The proof of Theorem \ref{OHBSC} relies on classical results in  \cite{Mok88}, \cite{Che07}, \cite{GZ10}, \cite{CT12} and \cite{FLW17} that a compact K\"ahler manifold with positive orthogonal holomorphic bisectional curvature must be biholomorphic to $\C\P^n$.    For more  discussions on compact K\"ahler manifolds with positive holomorphic sectional curvature, we refer to \cite{YZ19}, \cite{Yan20}, \cite{Yang2021}, \cite{Ni21}, \cite{Mat22},  \cite{NZ22}, \cite{LZZ21+}, \cite{ZZ23+} and the references therein.

\vskip 1\baselineskip

\noindent\textbf{Acknowledgements}. The second named author  would like to thank Bing-Long Chen, Jixiang Fu and Valentino Tosatti for
helpful discussions.  He would also like to thank   Professor Shing-Tung Yau and Professor Kefeng Liu for their support, encouragement and stimulating
discussions over many years. The second named author is partially supported by National Key R\&D
Program of China 2022YFA1005400 and NSFC grants (No. 12325103, No. 12171262
and No. 12141101). 

\vskip 2\baselineskip
	
	\section{Estimates of conjugate radius}
	
	In this section we obtain conjugate radius estimates for compact K\"ahler manifolds and establish Theorem \ref{main1}, Corollary \ref{main2} and Theorem \ref{main3}.	Let $(M,g)$ be a complete Riemannian manifold.  	For  each $t\in\R$, there is a flow induced by geodesics of $(M,g)$ \beq \phi_t:TM\to TM,\ \ 
	\phi_t(p, v)=\left(\gamma_v(t), \gamma_v'(t)\right) \label{realflow}
	\eeq 
	where $v\in T_p M$ and $\gamma_v(t)=\exp_p(tv)$.  We also write it as $\phi_t(v)=\gamma'_v(t)$ for simplicity. 
	We shall show that this flow is volume preserving, i.e.,  the determinant of the Jacobian map of $\phi_t$ with respect to the induced Sasaki metric on $TM$ is  $1$.
	Let's  describe the set up briefly and we refer to \cite[Lemma~3.1]{Gre63} and \cite{Gro16} for more details. 
	Let $\pi:TM\to M$ be the projection of the tangent bundle. There is a natural bundle map \beq C:TTM\to TM\eeq  which is defined as follows. Let $m$ be a point in $TM$.
	\bd 
	\item For $Z\in T_mTM$,  there exists a smooth curve $V:(-\eps,\eps)\to TM$ such that $V(0)=m$ and $V'(0)=Z$.
	\item Let $\gamma=\pi\circ V:(-\eps,\eps)\>M$ be a curve. The map $C$ is given by  \beq C(m, Z)=\left(\pi(m), \hat\nabla_{\left.\frac{d}{dt}\right|_{t=0}}\left(\gamma^*V\right)\right)\eeq  where $\hat\nabla$ is the pullback Levi-Civita connection along $\gamma:(-\eps,\eps)\>M$.
	\ed 
	It is easy to see that the bundle map $C$ is well-defined and smooth.
	Moreover, $$\cH:=\ker C \qtq{and} \cV:=\ker\pi_*$$ are subbundles of $TTM$ satisfying  $TTM=\cH\oplus \cV$.   It is well-known that there exists a unique Riemannian metric $\hat g$ on smooth manifold $TM$, which is called the \emph{Sasaki metric},  such that $\cH\perp\cV$ and  for all $p\in M$ and $v\in T_pM$, and the maps
	\beq 
	\pi_*:H_{(p,v)}\to T_pM\qtq{and}C:\cV_{(p,v)}\to T_pM
	\eeq 
	are linear isometries where $T_pM$ is endowed with the Euclidean metric induced by $g$.  Let $\hat g_p$ be the induced metric on the submanifold $T_pM$ of $(TM,\hat g)$. 
	\blemma \label{transition}  $\hat g_p$ coincides with the Euclidean metric on $T_pM$ induced by $g$. \elemma 
	\bproof 	Fix two points $v,w\in T_pM$ and set $$V(t)=v+tw:\R\>T_pM\subset TM.$$ Then $\gamma:=\pi\circ V:\R\>M$ is a constant, i.e. $\gamma(t)\equiv p$. Therefore $V'(0)\in T_{(p,v)}TM$ and it is in $$\ker\pi_{*(p,v)}=\cV_{(p,v)}.$$
	Since $V(t)$ is also a curve in $T_pM$,
	one can identify $V'(0)\in T_{v}(T_pM)\subset T_{(p,v)}TM$ and 
	\beq |V'(0)|_{\hat{g}_p}=|V'(0)|_{\hat{g}}.\eeq
	Let $\{x^i\}$ be local coordinates near $p\in M$, and $e_i=\frac{\p}{\p x^i}$ around $p$.
	We write $v=v^ie_i(p)$ and $w=w^je_j(p)$. 
	Then $\left(\gamma^* V\right)(t)=\left(v^i+tw^i\right) \hat{e_i}$ where $\hat{e_i}=\gamma^*e_i$ and 
	\beq 
	\hat\nabla_{\left.\frac{d}{dt}\right|_{t=0}}\left(\gamma^*V\right)
	=w^i\hat{e_i}(0)+v^i\hat\nabla_{\left.\frac{d}{dt}\right|_{t=0}}\hat{e_i}
	=w^i\hat{e_i}(0)=w
	\eeq 
	where we use the fact that $\gamma(t)\equiv p$.
	Therefore, 
	\beq 
	C(v,V'(0))=\left(p, \hat\nabla_{\left.\frac{d}{dt}\right|_{t=0}}\left(\gamma^*V\right)\right)
	=(p,w).
	\eeq 
	That is $C(V'(0))=w$.
	Since $C:\cV_{(p,v)}\to T_pM$ is a linear isometry, one has \beq |V'(0)|_{\hat{g}_p}=|V'(0)|_{\hat{g}}=|w|_g. \eeq 
	By using the identification $T_v(T_pM)\cong T_pM$,
	one deduces that the Riemannian metric $\hat{g}_p$ on $T_pM$ coincides with the Euclidean metric on $T_pM$ induced by $g$. 
	\eproof 
	
	\noindent	Let $(M,g)$ be a compact and oriented Riemannian manifold and   $SM$ be the unit tangent bundle of $M$.  For simplicity, the induced metric on the submanifold $SM$ of $TM$ is  denoted  by $\hat g$
	and the induced metric on  $S_pM$ is also denoted  by $\hat g_p$. By using Lemma \ref{transition}, one obtains	the following well-known lemma  (e.g. \cite{Gro16}) in Riemannian geometry.

	\blemma\label{prop geodesic flow}	For each $f\in C^\infty(SM, \R)$, one has 
	\beq 
	\int_{SM} f\,d\mathrm{vol}_{\hat g} =\int_M \left(\int_{S_pM} \left.f\right|_{S_pM}\,d\mathrm{vol}_{\hat g_p} \right)d\mathrm{vol}_g=\int_{SM}f \circ \phi_t \,d\mathrm{vol}_{\hat g}.
	\label{formula integration}
	\eeq 
	
	\elemma

	\vskip 1\baselineskip
	
	\noindent We introduce a complex analog of the flow (\ref{realflow}) on  a compact K\"ahler manifold $(M,\omega_g)$. For each $t\in \R$, there is an induced  flow on the holomorphic tangent bundle 	\beq  
	\psi_t=\Phi\circ\phi_t\circ\Phi^{-1}:T^{1,0}M\to T^{1,0}M
	\eeq 
	where  the identification $\Phi:T_\R M\to T^{1,0}M$ is given by  $\Phi(v) =\frac{1}{\sqrt{2}}\left(v-\sq Jv\right)$ and $\phi_t$ is defined in (\ref{realflow}). There is an induced Riemannian metric on smooth manifold $T^{1,0}M$ which is given by \beq \tilde g=(\Phi^{-1})^*\hat g \eeq   where $\hat g$ is the Sasaki metric on the real tangent bundle $T_\R M$ of $M$. Let $UM$ be the unit holomorphic tangent bundle of $T^{1,0}M$, $\tilde{g}$ be the Riemannian metric on $UM$ induced by $(T^{1,0}M,\tilde g)$, and $\tilde g_p$ be the Riemannian metric on the submanifold $U_pM\subset T_p^{1,0}M$ which coincides with the Euclidean metric on $T^{1,0}_pM$ induced by $g$ as shown in Lemma \ref{transition}.

	\bproposition\label{complex integration}
	For each $f\in C^\infty(UM, \R)$, one has 
	\beq 
	\int_{UM}f\,d\mathrm{vol}_{\tilde g}=\int_M\left( \int_{U_pM} \left.f\right|_{U_pM}\,d\mathrm{vol}_{\tilde g_p} \right)d\mathrm{vol}_g=\int_{UM}f \circ \psi_t \,d\mathrm{vol}_{\tilde g}.
	\label{formula complex geodesic flow}
	\eeq 
	\eproposition
	
	\bproof
	Since $\Phi$ is a smooth isometry which sends $SM$ to $UM$, by Lemma \ref{prop geodesic flow}, we have
	\be
	\int_{UM}f \circ \psi_t \,d\mathrm{vol}_{\tilde g}
	&=& \int_{SM} f \circ \psi_t \circ\Phi\,d\mathrm{vol}_{\hat g}
	=\int_{SM} f \circ \Phi\circ\phi_t \,d\mathrm{vol}_{\hat g}\\
	&=& \int_{SM} f \circ \Phi \,d\mathrm{vol}_{\hat g}
	= \int_{UM} f  \,d\mathrm{vol}_{\tilde g}.
	\ee
	Moreover, by  integration formula (\ref{formula integration}) , one has 
	\[ 
	\int_{UM} f  \,d\mathrm{vol}_{\tilde g}  
	=\int_{SM} f \circ \Phi \,d\mathrm{vol}_{\hat g}
	=\int_M \left(\int_{S_pM} \left.f\circ \Phi\right|_{S_pM}\,d\mathrm{vol}_{\hat g_p}\right) d\mathrm{vol}_g.
	\] 
	Note   that the restriction $\Phi:(S_pM,\hat g_p)\to (U_pM,\tilde g_p)$ is a smooth isometry, and so
	\[
	\int_{U_pM} \left.f\right|_{U_pM}\,d\mathrm{vol}_{\tilde g_p}
	=\int_{S_pM} \left.f\circ \Phi\right|_{S_pM}\,d\mathrm{vol}_{\hat g_p}.
	\] 
	Therefore we obtain the conclusion.
	\eproof

	\noindent Before giving the proof of Theorem \ref{main1}, we need some  algebraic calculations. 
	
	\blemma\label{lem average}
	Let $(\S^{2n-1},g_{\mathrm{can}})\subset\C^n$ be the round sphere. Then 
	\[ 
	\int_{\S^{n-1}} \xi^i\bar{\xi^j}\xi^k\bar{\xi^\ell}\,d\mathrm{vol}_{g_{\mathrm{can}}}=\frac{\mathrm{Vol}(\S^{2n-1})}{n(n+1)}\left(
	\delta_{ij}\delta_{k\ell}
	+\delta_{i\ell}\delta_{jk}
	\right)
	\] 
	where $\xi=(\xi^1,\cdots,\xi^n)\in \S^{2n-1}\subset \C^n$. 
	\elemma

	\blemma\label{prop scalar curvature}
	Let $(M^n,\om_g)$ be a K\"ahler manifold. Fix a point $p\in M$ and let $U_pM=\left\{ v\in T^{1,0}_p M:|v|_g=1\right\}$.
	Then 
	\beq 
	2s(p)=\frac{n(n+1)}{\mathrm{Vol}(\S^{2n-1})}\int_{U_pM} R(V,\bar V,V,\bar V)\,d\mathrm{vol}_{\tilde g_p}
	\eeq 
	where $\tilde g_p$  is the induced  metric on $T^{1,0}_pM$.
	\elemma
	\bproof 
	Let $\{e_i\}_{i=1}^{n}$ be an unitary basis of $T^{1,0}_pM$,
	and $R_{i\bar jk\bar \ell}:=R(e_i,\bar{e_j},e_k,\bar{e_\ell})$. If  $\xi=(\xi^1,\cdots,\xi^n)\in \S^{2n-1}\subset \C^n$ and $V=\xi^i e_i$, 
	then by Lemma \ref{lem average},
	\be 
	\int_{U_pM} R(V,\bar V,V,\bar V)\,d\mathrm{vol}_{\tilde g_p}
	&=&\int_{\S^{2n-1}} R_{i\bar jk\bar \ell} \xi^i\bar{\xi^j}\xi^k\bar{\xi^\ell}\,d\mathrm{vol}_{g_{\mathrm{can}}}\\
	&=&\frac{\mathrm{Vol}(\S^{2n-1})}{n(n+1)}\sum_{i,j,k,\ell=1}^{n}R_{i\bar jk\bar \ell} \left(
	\delta_{ij}\delta_{k\ell}
	+\delta_{i\ell}\delta_{jk}
	\right)\\
	&=& \frac{\mathrm{Vol}(\S^{2n-1})}{n(n+1)} 2s(p)
	\ee 
	where $s(p)$ is the scalar curvature of the K\"ahler metric at point $p\in M$. 
	\eproof

	\vskip 1\baselineskip

	\noindent \emph{Proof of Theorem \ref{main1}.}		   
	Suppose  that $\mathrm{conj}(M,g)=a<+\infty$.
	Let $\gamma:[0,a]\to M$ be an arbitrary unit speed geodesic with $\gamma(0)=p$ and $\gamma'(0)=v\in T_p M$. 
	Consider a normal variational vector field along $\gamma$ $$W(t)=\sin\left(\frac{\pi t}{a}\right)J\gamma'(t).$$ 
	Since $\mathrm{conj}(M,g)= a$, by  the index form theorem, one has  \beq I_\gamma(W,W)=\int_0^a\left\{\LL\hat{\nabla}_{\frac{d}{dt}}W,\hat{\nabla}_{\frac{d}{dt}}W\RL-R(W,\gamma',\gamma',W)\right\}dt\geq 0.\eeq This implies
	\beq 
	\int_0^a \sin^2\left(\frac{\pi t}{a}\right)R(J\gamma',\gamma',\gamma',J\gamma')\,dt
	\le\int_0^a \frac{\pi^2}{a^2}\cos^2\left(\frac{\pi t}{a}\right)\,dt=\frac{\pi^2}{2a}.
	\label{key estimate}
	\eeq
	By using the index form theorem again, one deduces that the identity holds if and only if $W(t)$ is a Jacobi field along $\gamma$.
	We write $V_t:=\frac{1}{\sqrt{2}}\left(\gamma'(t)-\sq J\gamma'(t)\right)$ for each $t\in [0,a]$, and set $V:=\frac{1}{\sqrt{2}}\left(v-\sq Jv\right)$.
	Then for each $t\in [0,a]$, one has $$\psi_t V=\Phi\circ\phi_t\circ \Phi^{-1}(V)=\Phi\circ \phi_t(v)=\Phi(\gamma'(t))=V_t.$$  
	On the other hand, a straightforward calculation shows
	\beq 
	R(J\gamma',\gamma',\gamma',J\gamma')(t)=R(V_t,\bar V_t, V_t,\bar V_t)=R(\psi_t V,\bar{\psi_tV},\psi_t V,\bar{\psi_tV}).
	\eeq 
	Therefore, (\ref{key estimate}) is equivalent to \beq \int_0^a \sin^2\left(\frac{\pi t}{a}\right)R(\psi_t V,\bar{\psi_tV},\psi_t V,\bar{\psi_tV})\,dt\le \frac{\pi^2}{2a}.\label{key estimate complex} \eeq 
	Since $p$ and $v$ are arbitrary, one deduces that (\ref{key estimate complex}) holds for all $p\in M$ and $V\in U_pM$. By using Proposition \ref{complex integration},
	 one can integrate (\ref{key estimate complex}) over $UM$ and obtain
	\be 
	\frac{\pi^2}{2a}\mathrm{Vol}(M,\omega_g)\mathrm{Vol}(\S^{2n-1})&=&\int_{UM}\frac{\pi^2}{2a}d\mathrm{vol}_{\tilde g}\\
	&\ge& \int_{UM}\left(\int_0^a \sin^2\left(\frac{\pi t}{a}\right)R(\psi_t V,\bar{\psi_tV},\psi_t V,\bar{\psi_tV})\,dt\right)d\mathrm{vol}_{\tilde g}\\
	&=& \int_0^a \sin^2\left(\frac{\pi t}{a}\right)\left(\int_{UM}R(\psi_t V,\bar{\psi_tV},\psi_t V,\bar{\psi_tV})\,d\mathrm{vol}_{\tilde g}\right)dt.
	%&=& \int_0^a \sin^2\left(\frac{\pi t}{a}\right)\,dt\cdot \frac{\mathrm{Vol}(\S^{2n-1})}{n(n+1)}\int_M \mathrm{scal}(p)\,d\mathrm{vol}_g(p)\\
	%&=& \frac{a}{2}\frac{\mathrm{Vol}(\S^{2n-1})}{n(n+1)}\int_M \mathrm{scal}(p)\,d\mathrm{vol}_g(p).
	\ee 
	Note that for each $t\in\R$, by Proposition \ref{complex integration} and Lemma \ref{prop scalar curvature}, one has
	\be 
	\int_{UM}R(\psi_t V,\bar{\psi_tV},\psi_t V,\bar{\psi_tV})\,d\mathrm{vol}_{\tilde g}
	&=& \int_{UM}R( V,\bar{V}, V,\bar{V})\,d\mathrm{vol}_{\tilde g}\\
	&=& \int_M \left(\int_{U_pM} R( V,\bar{V}, V,\bar{V})\,d\mathrm{vol}_{\tilde g_p}\right) d\mathrm{vol}_g\\
	&=& \frac{\mathrm{Vol}(\S^{2n-1})}{n(n+1)}\int_M 2s\,d\mathrm{vol}_g\\
	&=& \frac{4\pi \mathrm{Vol}(\S^{2n-1})}{(n+1)!}\int_M c_1(M)\wedge [\om_g]^{n-1}.
	\ee 
	Therefore,  one has
	\be 
	\frac{\pi^2}{2a}\mathrm{Vol}(M,\omega_g)\mathrm{Vol}(\S^{2n-1})
	%&\ge& \int_{UM}\int_0^a \sin^2\left(\frac{\pi t}{a}\right)R(\psi_t V,\bar{\psi_tV},\psi_t V,\bar{\psi_tV})\,dtd\mathrm{vol}_{\tilde g}(V)\\
	%&=& \int_0^a \sin^2\left(\frac{\pi t}{a}\right)\int_{UM}R(\psi_t V,\bar{\psi_tV},\psi_t V,\bar{\psi_tV})\,d\mathrm{vol}_{\tilde g}(V)dt\\
	&\ge & \int_0^a \sin^2\left(\frac{\pi t}{a}\right)\,dt\cdot \frac{4\pi \mathrm{Vol}(\S^{2n-1})}{(n+1)!}\int_M c_1(M)\wedge [\om_g]^{n-1}\\
	&=& \frac{2\pi a }{(n+1)!}\mathrm{Vol}(\S^{2n-1})\int_M c_1(M)\wedge [\om_g]^{n-1}.
	\ee 
	Thus we obtain the inequality (\ref{formula key conclusion}).  Furthermore,
	suppose that the identity in (\ref{formula key conclusion}) holds. One can deduce that the identity in (\ref{key estimate complex}) holds for all $p\in M$ and $V\in U_pM$. Moreover,  the identity in (\ref{key estimate}) holds for any unit-speed geodesic $\gamma:[0,a]\to M$, 
	and  $W(t)=\sin(\pi t/a)J\gamma'(t)$ is a Jacobi field along $\gamma$. This implies
	\be
	\left\langle\hat\nabla_{\frac{d}{dt}}\hat\nabla_{\frac{d}{dt}}W+R(W,\gamma')\gamma',W\right\rangle(t)
	=\sin^2\left(\frac{\pi t}{a}\right)\left(-\frac{\pi^2}{a^2}+R(J\gamma',\gamma',\gamma',J\gamma')(t)\right)=0.
	\ee 
	Therefore, for any $t\in (0,a)$,
	$$R(J\gamma',\gamma',\gamma',J\gamma')(t)=\frac{\pi^2}{a^2}$$
	and by continuity,  one obtains $R(Jv,v,v,Jv)=\pi^2/a^2$ where $v=\gamma'(0)\in T_p M$.
	Since $v$ and $p$ are arbitrary, we conclude that  $(M,\omega_g)$ has constant holomorphic sectional curvature $\pi^2/a^2$, and so $(M,\omega_g)$ is isometrically biholomorphic to $\left(\C\P^n,\frac{2a^2}{\pi^2}\omega_{\mathrm{FS}}\right)$.\\
	
	Suppose  that $\mathrm{conj}(M,g)=+\infty$. 	Let $\alpha$ be an arbitrary positive number,
	and $\gamma:[0,\alpha]\to M$ be a unit-speed geodesic. 
	Since $\mathrm{conj}(M,g)>\alpha$, by using the index form theorem, one has $I_\gamma(W,W)> 0$ where $W(t)=\sin(\pi t/\alpha)J\gamma'(t)$.
	By using similar arguments as above, one can show 
	\[
	\int_0^\alpha \sin^2\left(\frac{\pi t}{a}\right)R(J\gamma',\gamma',\gamma',J\gamma')\,dt< \frac{\pi^2}{2\alpha}.
	\]
	We can repeat previous arguments and obtain
	\[
	\frac{4\alpha^2}{\pi (n+1)!}\int_M c_1(M)\wedge [\omega_g]^{n-1} < \mathrm{Vol}(M,\omega_g).
	\]
	Since  $\alpha$ is arbitrary, we deduce that  \beq \int_M c_1(M)\wedge [\omega_g]^{n-1}\le 0.\eeq 
	Hence, the inequality in (\ref{formula key conclusion}) holds. Moreover, the identity in (\ref{formula key conclusion}) can not hold.
	\qed 
	
	\vskip 1\baselineskip
	\noindent \emph{Proof of Corollary \ref{main2}.}
	Since the scalar curvature $s\ge n(n+1)$, one has $$\int_M c_1(M)\wedge [\omega_g]^{n-1}\ge \frac{(n+1)!}{2\pi}\mathrm{Vol}(M,\om_g).$$
	Therefore, by Theorem \ref{main1}, one deduces that $a<+\infty$ and 
	\[
	\frac{2a^2}{\pi^2}\mathrm{Vol}(M,\om_g)\le \frac{4a^2}{\pi (n+1)!}\int_M c_1(M)\wedge [\omega_g]^{n-1} \le \mathrm{Vol}(M,\omega_g).
	\]
	This implies $a\le \pi/\sqrt{2}$.
	Moreover,  if $a=\pi/\sqrt{2}$, by the proof of Theorem \ref{main1}, one can see that $(M,\omega_g)$ is isometrically biholomorphic to $(\C\P^n,\om_{\mathrm{FS}})$.
	\qed

	\vskip 1 	\baselineskip
	
	\bproof[Proof of Theorem \ref{main3}]
	Suppose there exists some K\"ahler metric $\om_g$ on $M$ such that  $\mathrm{conj}(M,g)=+\infty$.
	By the proof of Theorem \ref{main1}, we deduce that $$\int_M c_1(M)\wedge [\omega_g]^{n-1} \le 0.$$ 
	By \cite[Theorem~1.1]{Yan19b}, we conclude that $K_M$ is pseudo-effective. However, by \cite[Theorem~1.5]{Yan19a}, if $K_M^{-1}$ is RC-positive, then $K_M$ is not pseudoeffecive and this is a contradiction.
	\eproof

\vskip 1\baselineskip
	
	\section{Injectivity radius, volume comparison and rigidity theorems for holomorphic sectional curvature}
	In this section, we investigate the geometry of complete K\"ahler manifolds with positive holomorphic sectional curvature (HSC) and demonstrate Theorem \ref{main4} and Theorem \ref{main5}.
	Let $(M,g)$ be a complete Riemannian manifold and $p\in M$. For small $r>0$, 
	$B_r(p):=\exp_p \left(B_r(0)\right)$
	is an open subset of $M$, and $\exp_p: B_r(0)\>B_r(p)$ is a diffeomorphism. The supremum of all such $r>0$ is called the {injectivity radius of $M$ at $p$} and it is denoted by $\mathrm{inj}_p(M,g)$. The {injectivity radius of $M$}, denoted by $\mathrm{inj}(M,g)$,  is $\inf_{p\in M} \mathrm{inj}_p(M,g)$.
	The following result is well-known  and we refer to \cite[pp.~274]{dC92} and \cite{Kli59}.
	\blemma\label{lem loop}
	Let $(M,g)$ be a complete Riemannian manifold and $p\in M$. Suppose that there exists some point $q\in\mathrm{cut}(p)$ such that $d(p,q)=d(p,\mathrm{cut}(p))=\ell$. Then one has
	\bd
	\item either $q$ is a conjugate point of $p$ along some minimizing geodesic from $p$ to $q$, 
	or there are exactly two unit-speed minimizing geodesics from $p$ to $q$, 
	say $\gamma_1,\gamma_2:[0,\ell]\to M$ such that $\gamma_1'(\ell)=-\gamma_2'(\ell)$;
	\item if in addition that $\mathrm{inj}_p(M,g)=\mathrm{inj}(M,g)$, and that $q$ is not conjugate to $p$ along any minimizing geodesic, 
	then there is a closed unit-speed geodesic $\gamma:[0,2\ell]\to M$ such that $\gamma(0)=\gamma(2\ell)=p$ and $\gamma(\ell)=q$.
	\ed
	\elemma
	
	\noindent We first show that on a compact K\"ahler manifold with positive holomorphic sectional curvature, the conjugate radius and the injectivity  radius are the same,  which is an analog of the classical result in \cite{Kli59} for  even dimensional orientable compact Riemannian manifolds with positive sectional curvature.
	\bproposition\label{prop inj=conj}
	Let $(M,\om_g)$ be a compact K\"ahler manifold with positive holomorphic sectional curvature.
	Then $\mathrm{inj}(M,g)=\mathrm{conj}(M,g)$.
	\eproposition
	
	\bproof
	Suppose for the sake of  contradiction that $\mathrm{inj}(M,g)<\mathrm{conj}(M,g)$.
	Since $M$ is compact, there exist $p\in M$ and $q\in \mathrm{cut}(p)$ such that $\ell:=d(p,q)=\mathrm{inj}(M,g)$.
	Since $\mathrm{inj}(M,g)<\mathrm{conj}(M,g)$, one deduces that $q$ is not conjugate to $p$ along any minimizing geodesic.
	Then by part $(2)$ of Lemma \ref{lem loop}, 
	there is a closed unit-speed geodesic $\gamma:[0,2\ell]\to M$ such that $$\gamma(0)=\gamma(2\ell)=p,\quad\gamma(\ell)=q.$$  In the following, we shall construct a third minimal geodesic connecting $q$ and $p$, and by part $(1)$ of Lemma \ref{lem loop}, this is a contradiction.\\

	Consider the variation $$\alpha:[0,1]\times [0,2\ell]\to M,\ \ \alpha(s,t)=\exp_{\gamma(t)}(s\cdot W(t)),$$
	where $W(t)=J\gamma'(t)$.
	Let $\bar\nabla$ and $\hat \nabla$ be the pullback Levi-Civita connections on $\alpha^*TM$ and $\gamma^*TM$ respectively. The first variation of the arclength of $\gamma(t)$ gives
	$$
	\left.\frac{d}{d s}\right|_{s=0} L\left(\alpha(s,\bullet)\right)
	=\left.\LL \gamma',W\RL\right|_{t=0}^{t=2\ell}-\int_0^{2\ell}\LL \hat\nabla_{\frac{d}{dt}}\gamma',W\RL dt
	=0.
	$$ 
	Since $\gamma'(0)=\gamma'(2\ell)$, $\alpha(s,0)=\alpha(s,2\ell)$ for all $s\in [0,1]$, and $\hat\nabla_{\frac{d}{dt}} J\gamma'=0$, the second variation of the arclength of $\gamma(t)$ is
	reduced to 
	\be
	\left.\frac{d^2}{d s^2}\right|_{s=0} L\left(\alpha(s,\bullet)\right)
	&=& \left.\LL \gamma',\left.\left(\bar\nabla_{\frac{\p}{\p s}}\alpha_*\left(\frac{\p}{\p s}\right)\right)\right|_{s=0}\RL\right|_{t=0}^{t=2\ell}
	-\int_0^{2\ell}R(J\gamma',\gamma',\gamma',J\gamma')\,dt\\&=& -\int_0^{2\ell}R(J\gamma',\gamma',\gamma',J\gamma')\,dt<0.
	\ee
 This implies that $\gamma$ is a local maximum of the arc-length functional.	We shall construct a minimal geodesic connecting $q$ and $p$. We write $\alpha_s(t)=\alpha(s,t)$, and it is clear that $$L(\alpha_s)<L(\gamma)=2\ell$$ for sufficiently small $s>0$.
	Let $p_s=\alpha_s(0)$, and $q_s=\alpha_s(t_s)$ be a point on the curve $\alpha_s$ that maximizes the distance to $p_s$. 
	By using this construction, one has $$ d(p_s,q_s)\le \frac{1}{2} L(\alpha_s)<\frac{1}{2}L(\gamma)=\ell=\mathrm{inj}(M,g)$$ for sufficiently small $s>0$.
	This implies that there exists a unique unit-speed minimal geodesic $\sigma_s:[0,\ell_s]\to M$ connecting $q_s$ and $p_s$.
		\begin{figure}[h]
		\centering
		\includegraphics[width=0.4\textwidth]{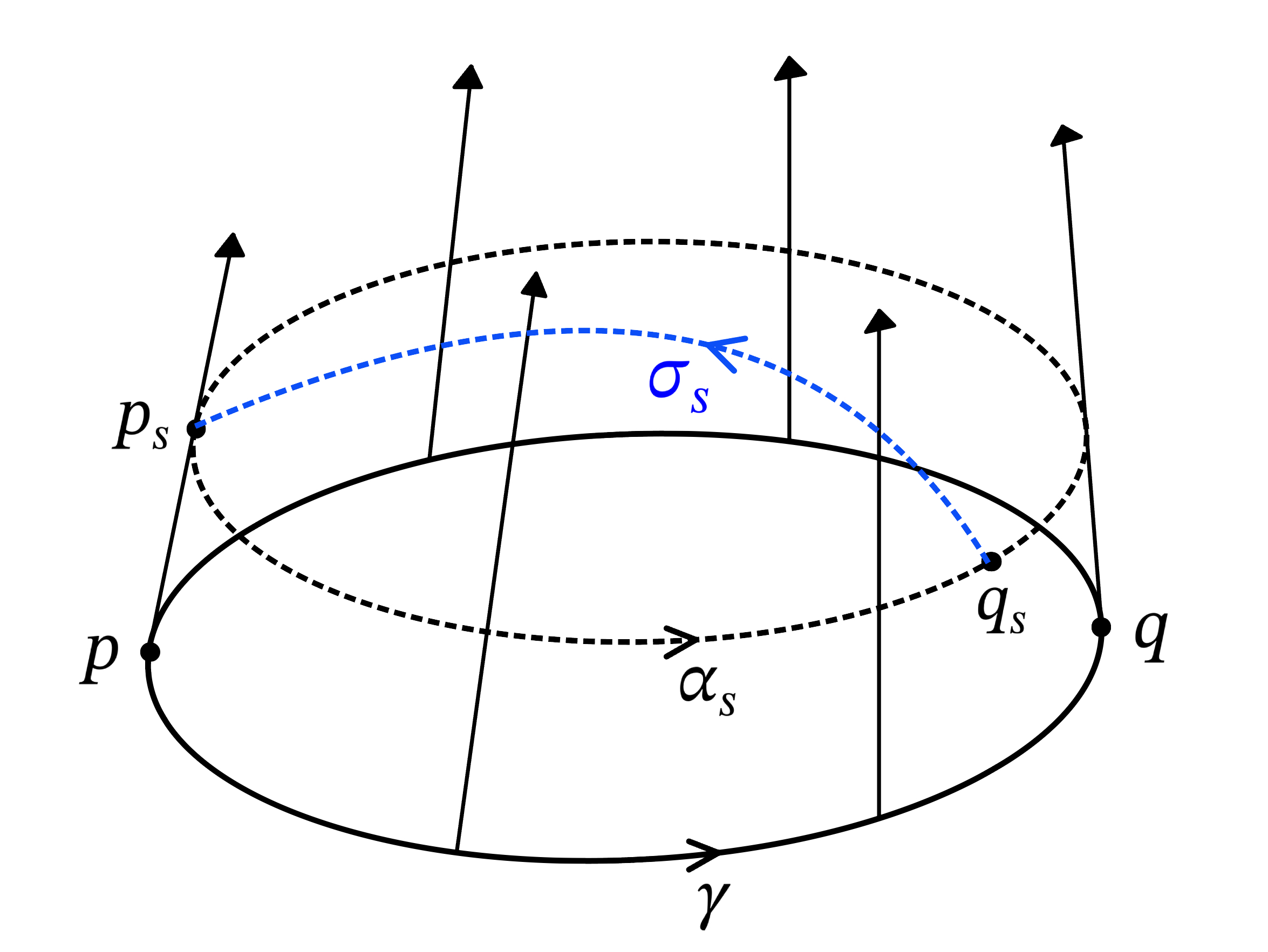}
		\caption{The variation  of $\sigma_s$}
	\end{figure}
	
\noindent	Moreover, there exists a smooth variation \beq \beta:(-\eps,\eps)\times [0,\ell_s]\to M\eeq  of $\sigma_s$ such that  for each $\tau\in(-\eps,\eps)$, the curve $\beta(\tau,\bullet)$ is a minimal geodesic with
	$$\beta(\tau,0)=\alpha_s(t_s+\tau)\qtq{and} \beta(\tau,\ell_s)=p_s.$$	Let $U(t)$ be the variational vector field of $\beta$. Then 
	$$
	U(0)=\left.\frac{d}{d\tau}\right|_{\tau=0}\alpha_s(t_s+\tau)=\alpha_s'(t_s) \qtq{and} U(\ell_s)=0.
	$$
	By the definition of $q_s$, one has 
	$$L(\sigma_s)=d(p_s,q_s)\geq d(p_s,\alpha_s(t_s+\tau))= L(\beta(\tau,\bullet))$$
	and so
	\beq 
	\left. \frac{d}{d\tau}\right|_{\tau=0}L\left(\beta(\tau,\bullet)\right)=0.
	\eeq

	\noindent On the other hand,  by the first variation formula,  for sufficiently small $s>0$,
	\beq 0=\left. \frac{d}{d\tau}\right|_{\tau=0}L\left(\beta(\tau,\bullet)\right)= \left.\LL\sigma_s',U\RL\right|_{t=0}^{t=\ell_s}= -\left\langle \sigma_s'(0),\alpha_s'(t_s)\right\rangle. \label{inj=conj 1} \eeq 
	Let $\{s_k\}$ be a sequence in the open interval $(0,1)$ which converges to $0$. 
	There exists a subsequence of $\{s_k\}$, which we also denote it by $\{s_k\}$,
	such that $$\lim_{k\to\infty}t_{s_k}=t_0$$ for some $t_0\in [0,2\ell]$. Thus, one has 
	$$
	\lim_{k\to\infty}q_{s_k}=\lim_{k\to\infty}\alpha(s_k,t_{s_k})=\alpha(0,t_0)=\gamma(t_0).
	$$
	Consider functions $f_k:[0,2\ell]\to \R$ and $f:[0,2\ell]\to \R$ given by $$f_k(t)=d(\alpha_{s_k}(t),p_{s_k})\qtq{and} f(t)=d(\gamma(t),p).$$
	One can see clearly that $f_k$ converges to $f$  uniformly, and 
	$$
	d(\gamma(t_0),p)
	=\lim_{k\to\infty}d(q_{s_k},p_{s_k})
	=\lim_{k\to\infty}\sup_{t\in[0,2\ell]}f_k(t)
	=\sup_{t\in[0,2\ell]}f(t)=\ell.
	$$
	Since $q$ is the only point on $\gamma$ that maximizes the distance to $p$, one deduces that
	\beq 
	\lim_{k\to\infty}t_{s_k}=t_0=\ell,\quad \lim_{k\to\infty}d(q_{s_k},p_{s_k})=\lim_{k\to\infty}\ell_{s_k}=\ell,\label{inj=conj 2}
	\eeq 
	and so $$\lim_{k\to\infty}\sigma_{s_k}(0)=\lim_{k\to\infty}q_{s_k}=\gamma(\ell)=q.$$
	Furthermore, by  compactness of the unit tangent bundle, 
	there exists a subsequence of $\{s_k\}$, which is also denoted by $\{s_k\}$, such that 
	\beq\lim_{k\to\infty}\sigma_{s_k}'(0)=w\label{inj=conj 3}\eeq
	for some $w\in T_qM$.
	We define a unit-speed geodesic \beq \sigma:[0,\ell]\to M,\ \ \ \sigma(t)=\exp_q(tw). \eeq 
	By continuity of the exponential map, one has
	\[ 
	\sigma(\ell)
	=\lim_{k\to\infty}\sigma_{s_k}(\ell)
	=\lim_{k\to\infty}\sigma_{s_k}(\ell_{s_k})
	=\lim_{k\to\infty} p_{s_k}=p.
	\] 
	Thus, $\sigma$ is also a minimal geodesic connecting $q$ and $p$.
	Moreover, by (\ref{inj=conj 1}), (\ref{inj=conj 2}) and (\ref{inj=conj 3}), one deduces that
	$$ 0=\lim_{k\to\infty} \LL \sigma_{s_k}'(0),\alpha_{s_k}'(t_{s_k})\RL= \LL \sigma'(0),\gamma'(\ell)\RL.$$
	Hence, $\sigma$ is a minimal geodesic connecting $q$ and $p$, which is different from two  minimal geodesics connecting $q$ and $p$ given by  $\gamma$. This is a contradiction.	\eproof

	\noindent\emph{Proof of Theorem \ref{main5}.}
	By \cite[Theorem~1]{Tsu57} and Proposition \ref{prop inj=conj}, we know $M$ is compact and  $$\mathrm{inj}(M,g)=\mathrm{conj}(M,g).$$
	On the other hand,  since $\mathrm{HSC}\ge 2$, by Lemma \ref{prop scalar curvature}, one deduces that $$s\ge n(n+1).$$
Now the estimate in \eqref{inj=conj} follows from Corollary \ref{main2},
	and the identity in \eqref{inj=conj} holds if and only if $(M,\omega_g)$ is isometrically biholomorphic to $(\C\P^n,\omega_{\mathrm{FS}})$.
	\qed

	\vskip 1\baselineskip
	\noindent Before proving Theorem \ref{main4}, we need the following result, which might be known to experts along this line. For the reader's convenience, we include a proof here.

	\bproposition\label{lem from Jacobi to curvature}
	Let $(M,\om_g)$ be a complete K\"ahler manifold, $p\in M$ and $U=M\backslash\mathrm{cut}{(p)}$.  Let  $\kappa\in \R$ and $\gamma:[0,\ell]\to M$ be an arbitrary unit-speed geodesic satisfying $\gamma(0)=p$ and $\gamma(t)\in U$ for all $t\in [0,\ell]$.
	Then the following statements are equivalent.
	\bd 
	\item  
	Every Jacobi field $J(t)$ along $\gamma$ with $J(0)=0$ and $\LL J,\gamma'\RL\equiv 0$ is of the form 
	\beq
	J(t)=a \operatorname{sn}_{\kappa/2}(t)E(t)+ b\operatorname{sn}_{2\kappa}(t)J\gamma'(t)\label{Jacobi field}
	\eeq
	where $E(t)$ is some parallel vector field along $\gamma$ with $\LL E(t),\gamma'(t)\RL=\LL  E(t),J\gamma'(t)\RL\equiv 0$ and $|E(t)|\equiv 1$.
	\item $(M,\om_g)$ has constant holomorphic bisectional curvature $\kappa$. %i.e. $R_{i\bar j k\bar l}=\kappa\left(g_{i\bar j}g_{k\bar l}+g_{i\bar l}g_{k\bar j}\right)$.
	\ed 
	\eproposition

	\bproof
	A straightforward calculation shows that $(2)$ implies $(1)$.  We shall show that (1) implies (2).
	Let $(N,\om_h)$ be a complete K\"ahler manifold with $\mathrm{HBSC}(N,\om_h)\equiv \kappa$.  Fix a point $\tilde p\in N$.
	In the following, we shall construct a holomorphic local isometry $\phi:(U,\omega_g)\to (N,\omega_h)$.
	This implies $\mathrm{HBSC}\equiv \kappa$ on $U$,
	and by continuity, we  conclude that $(M,\om_g)$ has constant holomorphic bisectional curvature $\kappa$. \\

	We choose a linear isometry $F:T_pM\to T_{\tilde p}N$ such that for each $v\in T_pM$,
	\beq F\left(J_Mv\right) =J_N\left(F(v)\right)\qtq{and} |F(v)|=|v|\label{F choice}\eeq
	where $J_M$ and $J_N$ are complex structures on $M$ and $N$ respectively.
	There is a smooth map given by
\beq 
		\phi= \tilde{\exp}_{\tilde p}\circ F\circ {\exp}_p^{-1}: U\to N.
\eeq 
	We claim that $\phi$ is a local isometry.
	Indeed, fix some $q\in U\setminus\{p\}$, and let $\gamma:[0,\ell]\to M$ be the unique unit-speed minimal geodesic connecting $p$ and $q$.
	We first show that for any $w\in T_qM$ with $\LL w,\gamma'(\ell)\RL=0$, one has 
	\beq |\phi_*w|=|w| . \label{formula phi local isometry}\eeq
  To this purpose,	let $J(t)$ be the unique Jacobi field along $\gamma$ with $J(0)=0$ and $J(\ell)=w$.
	Let $\tilde\gamma:[0,\ell]\to N$ be the geodesic given by $\tilde\gamma(t):=\phi(\gamma(t))$,
	and let $$ \tilde J(t):=\phi_*(J(t))$$ be a vector field along $\tilde{\gamma}$.
	One can see clearly that $\LL J,\gamma'\RL\equiv 0$ and $\tilde J(t)$ is a Jacobi field along $\tilde\gamma$ with $\tilde J(0)=0$ and $\tilde J'(0)=F(J'(0))$. 
	By using (\ref{F choice}), it is easy to see that $\LL\tilde J,\tilde\gamma'\RL\equiv 0$.
	%We shall show that $|\tilde J(t)|=|J(t)|$ for all $t\in [0,\ell]$.
	%The key point is that both Jacobi fields are of the form (\ref{Jacobi field}) with respect to their corresponding geodesics.
	Since $\mathrm{HBSC}(N,\om_h)\equiv \kappa$, by part $(1)$, the Jacobi field $\tilde J(t)$ is of the form
	\beq
	\tilde J(t)=\tilde a \operatorname{sn}_{\kappa/2}(t)\tilde E(t)+ \tilde b\operatorname{sn}_{2\kappa}(t)J_N\tilde \gamma'(t)\label{tilde J form}
	%\qtq{and}
	%|\tilde J(t)|^2={a^2 \operatorname{sn}_{\kappa/2}^2(t)+b^2\operatorname{sn}^2_{2\kappa}(t)}
	\eeq
	where $\tilde E(t)$ is some parallel vector field along $\tilde \gamma$ with $\LL \tilde E(t),\tilde \gamma'(t)\RL=\LL \tilde E(t),J_N\tilde \gamma'(t)\RL\equiv 0$ and $|\tilde E(t)|\equiv 1$.
	Moreover, by using (\ref{F choice}) and $\tilde J'(0)=F(J'(0))$, one deduces that
	\beq |J'(0)|=|\tilde J'(0)|\qtq{and}\LL J'(0),J_M\gamma'(0)\RL=\LL \tilde J'(0),J_N\tilde \gamma'(0)\RL.\label{Jacobi relation}\eeq
By assumption $(1)$,  $J(t)$ is a Jacobi field of the form (\ref{Jacobi field}), and a straightforward calculation shows that  $$a=\tilde a, \ \ \ b=\tilde b.$$ Hence,   $$|w|=|J(\ell)|=\sqrt{a^2 \operatorname{sn}_{\kappa/2}^2(\ell)+b^2\operatorname{sn}^2_{2\kappa}(\ell)}=|\tilde J(\ell)|=|\phi_*(J(\ell))|=|\phi_*w|.$$ 
	We complete the proof of (\ref{formula phi local isometry}). \\

\noindent	Moreover, for $w'\in T_qM$, it can be written as $w'=w+c\cdot\gamma'(\ell)$ where $\LL w,\gamma'(\ell)\RL=0$. Since $\LL \phi_*w,\tilde\gamma'(\ell)\RL=\LL\tilde J(\ell),\tilde\gamma'(\ell)\RL=0$,
one has
	$$ | \phi_* w'|^2
	=| \phi_*w+c\cdot\tilde\gamma'(\ell)|^2= |\phi_*w|^2+c^2=|w|^2+c^2=|w'|^2.$$
	Since $q$ and $w'$ are arbitrary, one obtains that $\phi$ is a local isometry.\\

	Furthermore, let $P_{\gamma}:T_pM\to T_qM$ and $P_{\tilde\gamma}:T_{\tilde p}N\to T_{\tilde\gamma(\ell)}N$ be parallel transports along $\gamma$ and $\tilde\gamma$ respectively.
	Since $\phi$ is a local isometry, for all $v\in T_pM$ one has
	\be 
		\phi_* \left(J_M\left(P_{\gamma}v\right)\right)
		&=&\phi_* \left(P_{\gamma}( J_Mv)\right)
		=P_{\tilde \gamma} \left(\phi_*( J_Mv)\right)
		=P_{\tilde \gamma} \left( J_N\left(\phi_*v\right)\right)\\
		&=&J_N\left(P_{\tilde \gamma} \left(\phi_*v\right)\right)
		=J_N\left(\phi_*\left(P_{\gamma} v\right)\right)
	\ee
	where we use (\ref{F choice}) and K\"ahler conditions that $\tilde{\nabla} J_N=0$ and $\nabla J_M=0$.
	Since $\gamma$ and $v$ are arbitrary, one gets
	$$\phi_*\left(J_Mw\right)= J_N\left(\phi_*(w)\right)$$ for all $w\in T_qM$ with $q\in U$.
	Therefore, $\phi$ is a holomorphic local isometry.
	\eproof

	\noindent The following lemma is the original  idea on RC-positivity which plays a key role in the proof of Theorem \ref{main4}, and we refer to \cite[Lemma~6.1]{Yan18} for more details.
	
	\blemma\label{lem variation}
	Let $(M,\omega_g)$ be a K\"ahler manifold and $p\in M$. 
	Let $e_1\in T^{1,0}_pM$ be a unit vector which minimizes the holomorphic sectional curvature $H$ of $\omega_g$ at point $p$.
	Then 
	\beq 
	2R(e_1,\bar e_1,W,\bar W)\geq \left(1+|\la W,e_1\ra|^2\right) R(e_1,\bar e_1,e_1,\bar e_1)
	\label{lem minimal vector}
	\eeq 
	for every unit vector $W\in T_p^{1,0}M$.
	\elemma

	\noindent \emph{Proof of Theorem \ref{main4}.}
	Let $\gamma:[0,\pi/\sqrt{2}]\to M$ be a unit-speed geodesic with $\gamma(0)=p$.
	Consider the normal variational vector field $$W(t)=\sin\left(\sqrt{2}t\right)J\gamma'(t).$$
	Since $\mathrm{conj}(M,p)\geq \pi/\sqrt{2}$, by the index form theorem, one has $$I_\gamma(W,W)\ge 0.$$
	On the other hand, by using $\mathrm{HSC}\ge 2$, one gets $R(J\gamma',\gamma',\gamma',J\gamma')\ge 2$, and this implies 
	\begin{align} 
	I_\gamma(W,W)
	&= \int_0^{\frac{\pi}{\sqrt{2}}} \left\{\LL\hat{\nabla}_{\frac{d}{dt}}W,\hat{\nabla}_{\frac{d}{dt}}W\RL-R(W,\gamma',\gamma', W)\right\}dt \nonumber\\
	&= \int_0^{\frac{\pi}{\sqrt{2}}}\left[ 2\cos^2\left(\sqrt{2}t\right)-\sin^2\left(\sqrt{2}t\right) \cdot R(J\gamma',\gamma',\gamma',J\gamma')(t)\right]dt  \nonumber\\
	&\le  \int_0^{\frac{\pi}{\sqrt{2}}} \left[2\cos^2\left(\sqrt{2}t\right)-2\sin^2\left(\sqrt{2}t\right) \right]dt=0.\label{index form}
	\end{align}
	Therefore, the identity in (\ref{index form}) holds, i.e. $I_\gamma(W,W)=0$.
	By the index form theorem again, one deduces that $W$ is a Jacobi field.
	The Jacobi field equation gives  \beq R(J\gamma',\gamma',\gamma',J\gamma')(t)\equiv 2 \label{curvature1} \eeq for  $t\in[0,\pi/\sqrt{2}]$.
	We write $V_t:=\frac{1}{\sqrt{2}}\left(\gamma'(t)-\sq J\gamma'(t)\right)\in T_{\gamma(t)}^{1,0}M$ for $t\in[0,\pi/\sqrt{2}]$. 
	It is clear that the holomorphic sectional curvature  $$ H(V_t)=R(V_t,\bar V_t,V_t,\bar V_t)\equiv 2.$$
	Thus one deduces that $V_t$ is a unit vector that minimizes the holomorphic sectional curvature $H$ at $\gamma(t)$. 
	By Lemma \ref{lem variation}, for any unit vector $U\in T_{\gamma(t)}^{1,0}M$, one has  \beq R(V_t,\bar V_t,U,\bar U)\ge 1+\left|\left\langle V_t,U\right\rangle\right|^2.\label{variation} \eeq 
	Let $\{e_1(t),\cdots,e_{2n}(t)\}$ be a parallel orthonormal frame along $\gamma$ such that $e_{2n}(t)=\gamma'(t)$ and for $1\leq k\leq n$ $$Je_{2k}(t)=e_{2k-1}(t) .$$ 
	If we set $U_k=\frac{1}{\sqrt{2}}\left(e_{2k}(t)-\sq e_{2k-1}(t)\right)$  for $1\leq k\leq n-1$,  by \eqref{variation}, one has
	\beq
	R(e_{2k-1},\gamma',\gamma',e_{2k-1})(t)+R(e_{2k},\gamma',\gamma',e_{2k})(t)=R(V_t,\bar V_t,U_k,\bar U_k)\ge 1.
	\label{curvature2}
	\eeq
	\vskip 1\baselineskip
\noindent 	Fix some $\ell\in (0,\pi/\sqrt{2})$.  Let $\sigma=\gamma|_{[0,\ell]}$. We define variational vector fields along $\sigma$:
	\beq \label{Jacobi fields for volume comparison}
	X_i(t)=\frac{\mathrm{sn}_{1/2}(t)}{\mathrm{sn}_{1/2}(\ell)}e_i(t),\quad1\le i\le 2n-2
	\qtq{and}X_{2n-1}(t)=\frac{\mathrm{sn}_{2}(t)}{\mathrm{sn}_{2}(\ell)}e_{2n-1}(t).
	\eeq 
	Let $\mathrm{cn}_k(t):=\mathrm{sn}_k'(t)$ and $\mathrm{ct}_k(t):=\mathrm{cn}_k(t)/\mathrm{sn}_k(t)$.
	By (\ref{curvature1}),
	\be
	I_\sigma(X_{2n-1},X_{2n-1})
	&=&  \int_0^\ell \left\{\LL\hat{\nabla}_{\frac{d}{dt}}X_{2n-1},\hat{\nabla}_{\frac{d}{dt}}X_{2n-1}\RL-R(X_{2n-1},\gamma',\gamma',X_{2n-1})\right\}dt\\
	&=&  \frac{1}{\mathrm{sn}_{2}^2(\ell)}\int_0^\ell \left[\mathrm{cn}_{2}(t)^2-\mathrm{sn}_{2}(t)^2 R(e_{2n-1},\gamma',\gamma',e_{2n-1})\right]dt\\
	&=& \frac{1}{\mathrm{sn}_{2}^2(\ell)} \int_0^\ell \left[\mathrm{cn}_{2}(t)^2- 2\mathrm{sn}_{2}(t)^2\right]dt 
	= \mathrm{ct}_{2}(\ell).
	\ee 
	By (\ref{curvature2}), for $1\le k\le n-1$,  one has
	\be
	\sum_{i=2k-1}^{2k}I_\sigma(X_i,X_i) 
	&=&  \sum_{i=2k-1}^{2k}\int_0^\ell \left\{\LL\hat{\nabla}_{\frac{d}{dt}}X_{i},\hat{\nabla}_{\frac{d}{dt}}X_{i}\RL-R(X_{i},\gamma',\gamma',X_{i})\right\}dt\\
	&=&  \frac{1}{\mathrm{sn}_{1/2}^2(\ell)}\int_0^\ell \left\{2\mathrm{cn}_{1/2}(t)^2- \mathrm{sn}_{1/2}(t)^2 \sum_{i=2k-1}^{2k} R(e_{i},\gamma',\gamma',e_{i})\right\}dt\\
	&\le& \frac{1}{\mathrm{sn}_{1/2}^2(\ell)} \int_0^\ell  \left[2\mathrm{cn}_{1/2}(t)^2- \mathrm{sn}_{1/2}(t)^2\right]dt 
	= 2\mathrm{ct}_{1/2}(\ell). 
	\ee
	Let $r(x)=d(p,x)$ be the distance function from point $p$. 	Suppose  $r$ is smooth at $\gamma(\ell)$. For $1\le i\le 2n-1$,
	let  $J_i$ be Jacobi fields along $\sigma$ such that $J_i(0)=0$ and $J_i(\ell)=e_i(\ell)$. It is well-known that (e.g. \cite[pp.~320--321]{Lee18})
	$$ \left(\Delta_g r\right)(\gamma(\ell))=  \sum_{i=1}^{2n-1}I_\sigma(J_{i},J_{i}). $$
	On the other hand, by the index form theorem, one has
	\begin{align}
	\left(\Delta_g r\right)(\gamma(\ell))
	&\le   I_\sigma(X_{2n-1},X_{2n-1})
	+ \sum_{k=1}^{n-1} \left[I_\sigma(X_{2k-1},X_{2k-1})+I_\sigma(X_{2k},X_{2k})\right]\label{volume -- index form theorem}\\
	&\le  \mathrm{ct}_{2}(\ell)+2(n-1)\mathrm{ct}_{1/2}(\ell).\nonumber
	\end{align}
	In the following we use similar arguments as in the proof of classical volume comparison theorems to make the conclusion. 
	Consider the map $$\Phi\colon \R^+\times\S^{2n-1}\to T_pM\setminus\{0\}\cong\R^{2n}\setminus\{0\}$$ given by $\Phi(t,v)=tv$, and define the volume density ratio as
	\beq 
	\lambda(t,v)=\frac{\chi_{\Sigma(p)}(tv)\cdot t^{2n-1}\sqrt{\det g}\circ\Phi(t,v)}{\mathrm{sn}^{2n-2}_{1/2}(t)\mathrm{sn}_{2}(t)},
	\eeq  	where $\Sigma(p)$ is the injectivity domain of $p$.
	Fix some $(\rho,\om)\in\R^+\times\S^{2n-1}$ and set $q:=\exp_p(\rho \om)$.
	If $\rho\omega\in \Sigma(p)$, then $r$ is smooth at $q$, and from the previous Laplacian estimate, one has 
	\begin{align}
	\left.\frac{\p}{\p t}\right|_{(\rho,\om)}\log \lambda
	&= \left.\frac{\p}{\p r}\right|_{q}\log\left(r^{2n-1}\sqrt{\det g}\right) - \left.\frac{d}{d t}\right|_{t=\rho}\log \left(\mathrm{sn}^{2n-2}_{1/2}(t)\mathrm{sn}_{2}(t)\right) \nonumber\\
	&= \left(\Delta r\right)(q) - \left(\mathrm{ct}_{2}(\rho)+2(n-1)\mathrm{ct}_{1/2}(\rho)\right) \le 0. \label{volume density}
	\end{align}
	If $\rho\omega\notin\Sigma(p)$, then $$\lambda(\rho,\om)=0.$$
	Thus, one deduces that for each $v\in \S^{2n-1}$, $\lambda(\bullet,v)$ is non-increasing for $t\in \R^+$. 
	Moreover, it is easy to see that for all $v\in \S^{2n-1}$, $\lim_{t\to 0^+}\lambda(t,v)=1.$
	Hence,  for any $(t,v)\in\R^+\times\S^{2n-1}$,
	\[ 
	\lambda(t,v)\le 1.
	\] 
	On the other hand, since $\mathrm{HSC}\ge 2$, one obtains $\mathrm{diam}(M,g)\le \pi/\sqrt{2}$,
	and so 
	\be
	\mathrm{Vol}(M,\om_g)
	&=& \int_{\S^{2n-1}}\int_0^{\frac{\pi}{\sqrt{2}}} \chi_{\Sigma(p)}t^{2n-1}\sqrt{\det g}\circ\Phi(t,v)\,dtd\mathrm{vol}_{\S^{2n-1}} \\
	&\le & \int_{\S^{2n-1}}\int_0^{\frac{\pi}{\sqrt{2}}} \mathrm{sn}^{2n-2}_{1/2}(t)\mathrm{sn}_{2}(t)\,dtd\mathrm{vol}_{\S^{2n-1}} 
	%= \mathrm{Vol}(\S^{2n-1})\int_0^{\frac{\pi}{\sqrt{2}}} \mathrm{sn}^{2n-2}_{1/2}(t)\mathrm{sn}_{2}(t)\, dt
	=\mathrm{Vol}(\C\P^n,\om_{\mathrm{FS}})
	\ee
	where $\Sigma(p)$ is the injectivity domain of $p$.
	This  establishes the inequality  (\ref{vol comparison}).  If the identity in (\ref{vol comparison}) holds,  it is clear that for all $ (t,v)\in (0,\pi/\sqrt{2})\times\S^{2n-1}$,  $$\lambda(t,v)\equiv 1.$$
	This implies that $$\Sigma(p)=B\left(0,\pi/\sqrt{2}\right),$$ and the identity in (\ref{volume density}) holds for all $q=\exp_p(\rho\om)$ with $(\rho,\om)\in (0,\pi/\sqrt{2})\times\S^{2n-1}.$
	In particular, if $\gamma\colon[0,\pi/\sqrt{2}]\to M$ is a unit-speed geodesic with $\gamma(0)=p$,
	and $\sigma=\gamma|_{[0,\ell]}$ for some $\ell\in (0,\pi/\sqrt{2})$,
	then the identity in (\ref{volume -- index form theorem}) holds.
	By the index form theorem, the vector fields $X_i$ given by (\ref{Jacobi fields for volume comparison}) are Jacobi fields along $\sigma$. We conclude that  every Jacobi field along $\sigma$ with $J(0)=0$ and $\LL J,\sigma'\RL\equiv 0$ is of the form
	\[
	J(t)=a \operatorname{sn}_{1/2}(t)E(t)+ b\operatorname{sn}_{2}(t)J\gamma'(t)
	\]
	where $E(t)$ is some parallel vector field along $\gamma$ with $\LL E(t),\gamma'(t)\RL=\LL  E(t),J\gamma'(t)\RL\equiv 0$ and $|E(t)|\equiv 1$.
	By Proposition \ref{lem from Jacobi to curvature}, one obtains that $(M,\om_g)$ has $\mathrm{HBSC}\equiv 1$,
	and so $(M,\om_g)$ is isometrically biholomorphic to $(\C\P^n,\om_{\mathrm{FS}})$. \qed

	\vskip 1\baselineskip 
	\noindent We propose the following problem for further investigation.
	
	\begin{problem}\label{HSC} Let $M$ be a compact K\"ahler manifold with positive holomorphic sectional curvature. Does the  volume comparison  theorem hold? How about the diameter and volume rigidity? 
	\end{problem}

	\vskip 1\baselineskip

	\section{Volume comparison and rigidity  theorems for orthogonal holomorphic bisectional curvature}

In this section, we investigate the geometry of  complete K\"ahler manifolds with positive orthogonal holomorphic bisectional curvature (OHBSC) and prove Theorem \ref{OHBSC}. Let $(M,\omega_g)$ be a complete K\"ahler manifold. Recall that, $(M,\omega_g)$ has $\mathrm{OHBSC}\geq 1$ if for any $p\in M$ and  unit vectors $X,Y\in T^{1,0}_pM$ with $g(X,\bar Y)=0$, one has 
\beq R(X,\bar X, Y,\bar Y)\geq 1.\eeq 
As an analog of Meyers' theorem, we show:
\blemma\label{lem OHBSC compact}
Let $(M,\om_g)$ be a complete K\"ahler manifold with $\dim_\C M\geq 2$.  If $(M,\om_g)$ has $\mathrm{OHBSC}\geq 1$, then $$\mathrm{diam}(M,g)\leq \sqrt 2\pi.$$
In particular,  $M$ is compact.
\elemma

\bproof  Suppose for the sake of contradiction  that there exist two points $p$ and $q$ with distance $d(p,q)=\ell>\sqrt 2 \pi$. 
Let $\gamma:[0,\ell]\to M$ be a unit-speed minimal geodesic such that $\gamma(0)=p$ and $\gamma(\ell)=q$. 
Let $E(t)$ be a parallel vector field along $\gamma$ such that 
$$\LL E(t),\gamma'(t)\RL=\LL E(t),J\gamma'(t)\RL=0\qtq{and} |E(t)|\equiv 1.$$
Consider two variation vector fields along $\gamma|_{[0,\sqrt{2}\pi]}$
$$V_1(t)= \operatorname{sn}_{1/2}(t)E(t)\qtq{and}V_2(t)= \operatorname{sn}_{1/2}(t)JE(t).$$
Let  $X_t,Y_t\in T^{1,0}_{\gamma(t)}M$ be unit vectors given by
\[ 
X_t=\frac{1}{\sqrt{2}}\left(\gamma'(t)-\sq J\gamma'(t)\right)
\qtq{and}
Y_t=\frac{1}{\sqrt{2}}\left(E(t)-\sq JE(t)\right).
\]
Since $\mathrm{OHBSC}\ge 1$ and $g(X_t,\bar Y_t)=0$, a straightforward calculation shows that 
\[ 
R(E(t),\gamma'(t),\gamma'(t),E(t))
+R(JE(t),\gamma'(t),\gamma'(t),JE(t))=R(X_t,\bar X_t,Y_t,\bar Y_t)\ge 1.
\] 
Therefore, one obtains
\be 
\sum_{i=1}^{2} I_{\gamma|_{[0,\sqrt 2\pi]}}(V_{i},V_{i})
&=& \sum_{i=1}^{2}\int_0^{\sqrt{2}\pi} \left\{\LL \hat\nabla_{\frac{d}{dt}}V_{i}(t),\hat\nabla_{\frac{d}{dt}}V_{i}(t)\RL-R(V_{i}(t),\gamma'(t),\gamma'(t),V_{i}(t))\right\} dt\\
&=& \int_0^{\sqrt{2}\pi}\left[ 2\operatorname{cn}_{1/2}^2(t)-\operatorname{sn}_{1/2}^2(t)R(X_t,\bar X_t,Y_t,\bar Y_t)\right]dt\\
&\le& \int_0^{\sqrt{2}\pi} \left[2\operatorname{cn}_{1/2}^2(t) -\operatorname{sn}_{1/2}^2(t) \right]dt
= 0.
\ee 
By the index form theorem,  along the curve $\gamma|_{[0,\sqrt{2}\pi]}$, $\gamma(0)$ has a conjugate point $\gamma(t_0)$ for some $t_0\in (0,\sqrt 2 \pi]$. In particular,  $\gamma:[0,\ell]\>M$ is not a minimal geodesic,
and this is a contradiction.
Hence, we deduce that $$\mathrm{diam}(M,g)\leq \sqrt{2}\pi$$ and in particular $M$ is compact.
\eproof

\bremark It is well-known ( e.g. \cite{Mok88}, \cite{Che07}, \cite{GZ10}, \cite{CT12} and \cite{FLW17}) that  a compact K\"ahler manifold with positive orthogonal holomorphic bisectional curvature is biholomorphic to $\C\P^n$. We know from Lemma \ref{lem OHBSC compact} that $M$ is actually biholomorphic to $\C\P^n$,  and so  the diameter upper bound $\sqrt 2\pi$ is not sharp. 
\eremark

\noindent 
The following result is essentially known in some special cases (e.g. \cite{GZ10}, \cite{CT12}, \cite{FLW17}, \cite{NZ18}) and we present a proof here for the sake of completeness.
\blemma \label{prop OHBSC scalar}
Let $(M,\om_g)$ be a K\"ahler manifold with $\dim_\C M=n\ge 2$. If there exist two constants $a$ and $b$ such that $a\le \mathrm{OHBSC}\le b$, then  the scalar curvature $s$ satisfies
\[ 
n(n+1)a\le s\le n(n+1)b.
\] 
\elemma

\bproof
Suppose that $\left\{e_\alpha\right\}$ is an orthonormal basis of $T_p^{1,0} M$. 
Then one has
$$
R\left(e_\alpha-e_\beta, \overline{e_\alpha}-\overline{e_\beta}, e_\alpha+e_\beta, \overline{e_\alpha}+\overline{e_\beta}\right)
=R_{\alpha \bar{\alpha} \alpha \bar{\alpha}}+R_{\beta \bar{\beta} \beta \bar{\beta}}
-R_{\alpha \bar{\beta} \alpha \bar{\beta}}-R_{\beta \bar{\alpha} \beta \bar{\alpha}}
\ge 4 a
$$
for any $\alpha \neq \beta$.
Similarly, replacing $e_\beta$ by $\sqrt{-1} e_\beta$, one gets
$$
R_{\alpha \bar{\alpha} \alpha \bar{\alpha}}+R_{\beta \bar{\beta} \beta \bar{\beta}}
+R_{\alpha \bar{\beta} \alpha \bar{\beta}}+R_{\beta \bar{\alpha} \beta \bar{\alpha}}
\ge 4 a 
$$	for any $\alpha \neq \beta$.
The summation of  two inequalities gives
$$
R_{\alpha \bar{\alpha} \alpha \bar{\alpha}}+R_{\beta \bar{\beta} \beta \bar{\beta}}\ge 4a 
$$
for any $\alpha \neq \beta$.	This implies that
\be
s(p) & =& \sum_{\alpha, \beta} R_{\alpha \bar{\alpha} \beta \bar{\beta}}=\sum_\alpha \sum_{\beta \neq \alpha} R_{\alpha \bar{\alpha} \beta \bar{\beta}}+\sum_\alpha R_{\alpha \bar{\alpha} \alpha \bar{\alpha}} \\
& =& \sum_\alpha \sum_{\beta \neq \alpha} R_{\alpha \bar{\alpha} \beta \bar{\beta}}
+\frac{1}{2}\left(\sum_{\alpha=1}^n R_{\alpha \bar{\alpha} \alpha \bar{\alpha}}
+\sum_{\beta=1}^n R_{\beta \bar{\beta} \beta \bar{\beta}}\right) \\
& =&\sum_\alpha \sum_{\beta \neq \alpha} R_{\alpha \bar{\alpha} \beta \bar{\beta}}
+\frac{1}{2} \sum_{\gamma=1}^{n-1}\left(R_{\gamma \bar{\gamma} \gamma \bar{\gamma}}
+R_{\gamma+1 \overline{\gamma+1} \gamma+1 \overline{\gamma+1}}\right)
+\frac{1}{2}\left(R_{n \bar{n} n \bar{n}}+R_{1 \overline{1} 1 \overline{1}}\right) \\
& \ge&n(n-1)a+2na= n(n+1)a .
\ee
Hence, $s\ge n(n+1)a$.  The proof of the other part is similar. 
\eproof

\noindent 	\emph{Proof of Theorem \ref{OHBSC}.}		By Lemma \ref{lem OHBSC compact}, $M$ is compact. Moreover, since $\mathrm{OHBSC}>0$,
by \cite[Theorem~3.2]{GZ10},  one deduces that $M$ is biholomorphic to $\C\P^n$.
In particular, $H_{\bp}^{1,1}(M,\R)=\R$ and it is well-known that 
\[  
c^n_1(M)
= \int_{M}\left(\frac{\mathrm{Ric}(\om_{\mathrm{FS}})}{2\pi}\right)^n
= \left(\frac{n+1}{2\pi} \right)^n\int_{M}\om_{\mathrm{FS}}^n
= n!\left(\frac{n+1}{2\pi} \right)^n \mathrm{Vol}(\C\P^n,\om_{\mathrm{FS}}).
\]
Since $c_1(M) \in H_{\bp}^{1,1}(M, \mathbb{R})$,  one has $c_1(M)=\lambda [\omega_g]$ for some $\lambda\in\R$ and
\beq
\lambda \int_M \om_g^n=\int_M c_1(M) \wedge \om_g^{n-1}=\frac{1}{2\pi n} \int_M s\,\om_g^n.
\eeq
By Lemma \ref{prop OHBSC scalar}, one has $s\ge n(n+1)$, and so 
\[
\lambda =\frac{1}{2\pi n} \frac{\int_M s\,\om_g^n}{\int_M \om_g^n}\ge \frac{n+1}{2\pi}.
\]
This implies
$$ \mathrm{Vol}(M,\om_g)
=\frac{1}{n!}\int_M\om_g^n
=\frac{c^n_1(M)}{n!\lambda^n}
\le \left(\frac{2\pi }{n+1} \right)^n \frac{c^n_1(M)}{n!}=\mathrm{Vol}(\C\P^n,\om_{\mathrm{FS}}).$$
This is (\ref{formula HBSC vol upper}).
Suppose the identity in (\ref{formula HBSC vol upper}) holds. It is clear that 
$$\lambda=\frac{n+1}{2\pi}$$
and $s\equiv n(n+1)$.	It follows that 
$$
\left[\mathrm{Ric}(\om_g)\right]=2\pi c_1(M)=2\pi\lambda[\om_g]=(n+1)[\om_g].
$$
By $\p\bp$-lemma, there exists some $\phi\in C^\infty(M,\R)$ such that 
$$
\mathrm{Ric}(\om_g)=(n+1)\om_g+\sq\p\bp\phi.
$$
By taking trace, one deduces  that 
$$
\mathrm{tr}_{\om_g}\sq\p\bp\phi=s-n(n+1)\equiv 0.
$$
In particular, $\phi$ is a constant, and so $$\mathrm{Ric}(\om_g)=(n+1)\om_g.$$
By uniqueness of K\"ahler-Einstein metrics on $\C\P^n$, one obtains $\om_g=\Phi^*\om_{\mathrm{FS}}$ for some $\Phi\in \mathrm{Aut}(\P^n)$. 
Therefore, the the identity in (\ref{formula HBSC vol upper}) holds if and only if $(M,\om_g)$ is isometrically biholomorphic to $(\C\P^n,\om_{\mathrm{FS}})$. \qed 

\vskip 1\baselineskip

\noindent By using similar arguments, we also obtain the following volume comparison and rigidity result for complete K\"ahler manifolds with pinched orthogonal holomorphic bisectional curvature.
\btheorem	Let $(M,\om_g)$ be a complete K\"ahler manifold with dimension $n\ge 2$. If  $1\leq \mathrm{OHBSC}\leq a$ for some constant $a\geq 1$, then $M$ is compact and
\beq 
\mathrm{Vol}(\C\P^n,a^{-1}\om_{\mathrm{FS}})\leq 	\mathrm{Vol}(M,\om_g)\le  \mathrm{Vol}(\C\P^n,\om_{\mathrm{FS}}),
\eeq
and the first identity holds if and only if $(M,\om_g)$ is isometrically biholomorphic to  $(\C\P^n,a^{-1}\om_{\mathrm{FS}})$.
\etheorem

	\noindent
	\def\cprime{$'$} %newcommand of prime over letters

\end{document}